\documentclass[12pt,a4paper,reqno]{article}
\usepackage{amsmath}
\usepackage{amssymb}
\usepackage{amsthm}

\renewcommand{\(}{\left(}
\renewcommand{\)}{\right)}
\newcommand{\g}{\mathfrak{g}}
\newcommand{\hh}{\mathfrak{h}}
\newcommand{\ssll}{\mathfrak{sl}}
\newcommand{\Lie}{\mathfrak}
\newcommand{\K}{\mathbb{K}}
\newcommand{\U}{\overline{U}_q (\g)}
\newcommand{\Utwo}{\overline{U}_q(\ssll (2))}
\newcommand{\Uthree}{\overline{U}_q(\ssll (3))}
\newcommand{\UU}{\mathcal{U}}
\newcommand{\A}{\mathcal{A}}
\newcommand{\B}{\mathcal{B}}
\newcommand{\C}{\mathcal{C}}
\newcommand{\F}{\mathcal{F}}
\renewcommand{\L}{\mathcal{L}}

\newcommand{\ox}{\otimes}
\newcommand{\<}{\langle}
\renewcommand{\>}{\rangle}
\newcommand{\half}{\frac{1}{2}}

\newcommand{\third}{\frac{1}{3}}
\newcommand{\thrd}[1]{\frac{#1}{3}}
\newcommand{\qqquad}{\quad\quad\quad}
\newcommand{\qqqquad}{\quad\quad\quad\quad}
\newcommand{\imp}{\Rightarrow}
\newcommand{\impl}{\Longrightarrow}
\newcommand{\ad}{\operatorname{ad}}
\newcommand{\id}{\operatorname{id}}
\newcommand{\rad}{\operatorname{rad}}
\newcommand{\End}{\operatorname{End}}
\theoremstyle{plain} %% This is the default
\newtheorem{thm}{Theorem}[section]

\newtheorem{lemma}{Lemma}[section]
\newtheorem{propn}{Proposition}

\theoremstyle{definition}
\newtheorem{defn}{Definition}[section]

\theoremstyle{remark}

\numberwithin{equation}{section}

\newcommand{\thmref}[1]{Theorem~\ref{#1}}

\newcommand{\lemref}[1]{Lemma~\ref{#1}}

\begin{document}
\title{\begin{flushright}\small q-alg/9510004 \end{flushright}
  \smallskip
  Quantum Lie Algebras of Type $A_n$}
\author{Volodimir Lyubashenko\thanks{Current address: Mathematical
    Methods of Systems Analysis, Department of Applied Mathematics, Kiev
    Polytechnic Institute, prosp. Peremogy 37, 252056 Kiev, Ukraine}
    \and Anthony Sudbery\\[10pt] \small Department of Mathematics,
    University of York,\\[-2pt] \small Heslington, York, England YO1 5DD\\
\small  E-mail:  as2@york.ac.uk}
\date{10 October 1995}

\maketitle
\bigskip
\begin{abstract}
It is shown that the quantised enveloping algebra of $\ssll (n)$ contains
a quantum Lie algebra, defined by means of axioms similar to Woronowicz's.
This gives rise to Lie algebra-like generators and relations for the locally
finite part of the quantised enveloping algebra, and suggests a canonical
Poincar\'e-Birkhoff-Witt basis.
\end{abstract}
\bigskip
\section{Introduction}
It is generally believed that quantised enveloping algebras are not enveloping

algebras; Chari and Pressley, for example, assert \cite[p. 258]{ChaP}
that ``there
is no quantum Lie algebra'' underlying a quantised enveloping algebra
$U_q(\g)$. The purpose of this paper is to challenge this belief.
The ultimate intention is to describe all quantised enveloping algebras
as genuine enveloping algebras --- that is, as associative algebras generated
by a finite-dimensional system like a Lie algebra, with relations obtained
by interpreting the Lie brackets as quadratic expressions like commutators.
Here we succeed, for
the case where $\g$ is of type $A_n$, in finding generators and relations
which are essentially of the desired form; moreover, they appear to yield
a Poincar\'e-Birkhoff-Witt theorem for the enveloping algebra. They are
based on Lie brackets defined by means of the adjoint action of a Hopf
algebra on itself, the quantum Lie algebra being a certain ad-invariant
subspace of the quantised enveloping algebra. A generalisation of the
classical concept of Lie algebra emerges, which is close to that
found by Woronowicz \cite{Wor} in his theory of non-commutative differential
geometry on quantum groups. However, Woronowicz's theory yields a satisfactory
deformation of the classical Lie algebra only in the case of general linear
groups \cite{Diffform, SchWatZum, SS1}; in all other cases (see
\cite{Jurco, Karpacz,
CSSW, SS2}) the quantum
Lie algebra has a different dimension from the corresponding classical
one: it is always $n^2$ where $n$ is the dimension of a particular
representation.
In particular, quantum Lie algebras corresponding to $\ssll (n)$ and having
dimension $n^2 - 1$, such as are described in this paper, appear to be
new.

The following reservations must be made:

1. The relations for the quantised enveloping algebra are not, as expected,
quadratic-linear (deformed commutator = Lie algebra element) but homogeneous
quadratic, the right-hand side containing an extra central element as
a factor. To put it another way, the structure constants are multiples
of a
function of the Casimirs. This is an interesting feature; it means
that the PBW theorem leads to a description of the vector space structure
of the enveloping algebra as a space of polynomial functions not on a
flat space (namely the Lie algebra) but on a hypersurface of the same
dimension. Thus quantisation gives rise to curvature of the enveloping
algebra.

2. The algebra generated by our quantum Lie algebra cannot be expected
to be the same size as the quantised enveloping algebra, since this differs
from the classical enveloping algebra in containing infinite power series
(specifically, exponentials) in the generators. On the other hand, if,
as is common in the mathematical literature, one takes the generators
of the quantised enveloping algebra to be $q^{H_i}$ in place of the Cartan
subalgebra elements $H_i$, then one loses the polynomials in the $H_i$.
It is a remarkable fact \cite{JL} that the quantised enveloping algebra
so defined contains a subalgebra, its locally finite part, which has the
same structure as the classical enveloping algebra
and is a large part of the quantised enveloping algebra. It seems likely
that the enveloping algebra of the quantum Lie algebra, as defined in
this paper, coincides with the locally finite part of the quantised enveloping
algebra. In this paper this is proved for the case of $\Lie{sl}(2)$.

3. The construction of a quantum Lie algebra given here makes sense for
any simple Lie algebra, but it is only for the Lie algebras of type $A_n$
that it yields a system with the same dimension as the classical Lie algebra.
Thus it is only in these cases that a PBW theorem is likely to hold. Even
here, I have no general proof of the theorem; it can be checked by explicit
calculation for $A_1$ and $A_2$.

4. As indicated above, the commutator-like relations obtained from the
quantum Lie algebra are not the only relations in the quantised enveloping
algebra; there is also a relation between central elements, giving the
extra (non-classical) central generator in terms of more familiar Casimir
elements which are polynomials in Lie algebra elements. The existence
of this extra relation is at present only a conjecture in all cases except
$\Lie{sl}(2)$, for which it can be explicitly demonstrated.

The organisation of the paper is as follows. The next section contains
the definition of a quantum Lie algebra and the proof that such a system
can be found inside the quantised enveloping algebra. This relies heavily
on the theorems of Joseph and Letzter \cite{JL} concerning the locally
finite part of a quantised enveloping algebra. In section 3 we complete
the theory for $\Lie{sl}(2)_q$, proving the Poincar\'{e}-Birkhoff-Witt theorem,
determining the extra relation between central elements, and proving that
the algebra defined by the quantum Lie algebra relations, together with
the extra relation, is the locally finite part of the quantised enveloping
algebra. Section 4 contains some material on the quantum Lie algebra
$\Lie{sl}(3)_q$,
with formulae for the Lie brackets and a discussion of the rather surprising
symmetry properties of the enveloping algebra.

Delius et al. \cite{DelHuff, Del2} have also looked for quantum deformations
of Lie algebras in the ad-invariant subspaces of quantised enveloping
algebras. Their quantum Lie algebras lack the antisymmetry and Jacobi
axioms that are given here; on the other hand, they retain more of the
symmetry of classical Lie algebras. Majid \cite{Maj:Liealg} has proposed
axioms for ``braided Lie algebras''; these include a Jacobi identity but
no anticommutativity
axiom and, not being satisfied in the classical case, do not
constitute a generalisation of the Lie algebra axioms.

I am grateful to Volodimir Lyubashenko for enlightening discussions.

{\em Notation.}\quad I use the usual notation for coproducts in a Hopf
algebra:
\begin{eqnarray*}
\Delta (x)&=&\sum x_{(1)}\ox x_{(2)},\\
(\Delta \ox \id)\circ\Delta (x)&=&\sum x_{(1)}\ox x_{(2)}\ox x_{(3)}\qquad
\text{etc.,}
\end{eqnarray*}
and for the $q$-analogue function:
\[ [x]_q=\frac{q^x-q^{-x}}{q-q^{-1}}. \]

\section{General Construction of Quantum Lie Algebras}
\begin{defn} A {\em quantum Lie algebra} over a field $\K$ is a $\K$-vector
space $L$ together with linear maps $\beta : L\ox L \to L$ (the quantum
Lie bracket, written $x \ox y \mapsto [x,y]$) and $\gamma : L \ox L \to
L \ox L$ (the quantum antisymmetriser) satisfying:

\noindent 1. Quantum antisymmetry: For $t \in L \ox L$,
\begin{equation}
\gamma(t)=0\; \imp\; \beta (t)=0.
\end{equation}

\noindent 2. The quantum Jacobi identity: For $x \in L$, define $\ad x:
L \to L$ by $\ad x(y)=[x,y]$. Then
\begin{equation} \label{antisym}
\ad [x,y]=m\circ (\ad \ox \ad )\circ \gamma (x \ox y)
\end{equation}
where $m : \End L \ox \End L \to \End L$ denotes multiplication of linear
maps on $L$.
\end{defn}
\begin{defn} A quantum Lie algebra is {\em balanced} if it satisfies a
second Jacobi
identity,

\noindent 3. The right quantum Jacobi identity: For $x \in L$, define $\rad x:
L \to L$ by $\rad x(y)=[y,x]$. Then
\begin{equation} \label{Jacobi}
\rad [x,y]=\,\stackrel{\leftarrow}{m}\circ (\rad \ox \rad )\circ \gamma
(x \ox y)
\end{equation}
where $\stackrel{\leftarrow}{m}:\End L \ox \End L \to \End L$ denotes
the opposite multiplication of linear maps: $\stackrel{\leftarrow}{m}(A\ox
B)=BA$.
\end{defn}

If $\gamma$ is the usual antisymmetrisation map, $\gamma (x \ox y)=x\ox
y-y\ox x$, the above become the usual antisymmetry and Jacobi axioms for
a Lie algebra. These generalised axioms were found by Woronowicz to hold
for the geometrical Lie bracket of left-invariant vector fields on a quantum
Lie group \cite{Wor}; he also found that in this case the quantum
antisymmetriser
was of the form $\gamma = 1-\sigma$ where $\sigma$ satisfies the braid
relation. This is not true of the quantum Lie algebras described in this
paper; I do not know if there is a weaker axiom of this type that they
satisfy.

Woronowicz included only one Jacobi identity among his axioms; the geometrical
construction does not obviously yield a balanced quantum Lie algebra.
This is also, at first sight, true of the construction from universal
enveloping algebras to be described in this paper. However, inspection
of the quantum Lie algebra $\ssll (2)_q$ reveals that it is in fact balanced,
and it seems likely that this is true more generally.

\begin{defn} A {\em representation} of a quantum Lie algebra $(L,
\beta, \gamma)$ is a linear map $\rho : L \to \End V$ for some vector space
$V$,
such that
\begin{equation}
\rho \circ \beta = m_V\circ (\rho \ox \rho)\circ \gamma : L \ox L \to
\End V
\end{equation}
where $m_V$ denotes multiplication of linear operators on $V$.
\end{defn}

Then the Jacobi identity \eqref{Jacobi} states that ad is a representation
of $L$ in $\End_{\K}V$, as with classical Lie algebras.

\begin{defn} The {\em universal enveloping algebra} of a quantum
Lie algebra $(L, \beta , \gamma)$ is the quotient of the tensor algebra
of $L$ by the ideal generated by $(\beta - \gamma)(V \ox V)$. Then every
representation of $L$ automatically generates a representation of its
universal enveloping algebra.
\end{defn}

Let $x_1, \ldots ,x_n$ be a basis for $L$, and let $\gamma _{ij}^{kl}$
be the matrix elements of $\gamma$ with respect to the basis ${x_i \ox
x_j}$ of $L \ox L$; in the classical case,
\[ \gamma _{ij}^{kl}=\delta_i^k\delta_j^l - \delta_i^l/\delta_j^k. \]
Then the Jacobi identity can be written as
\[ [\, [x_i, x_j], x_k]=\gamma_{ij}^{lm}[x_l, [x_m, x_k]\,], \]
the right Jacobi identity as
\[ [x_i, [x_j, x_k]\, ]=\gamma_{jk}^{lm}[\, [x_i, x_l], x_m],
\]
the representation property as
\[ \rho ([x_i,x_j])=\gamma_{ij}^{kl}\rho(x_k)\rho(x_l) \]
and the relations in the universal enveloping algebra as
\[ \gamma_{ij}^{kl}x_k x_l = [x_i, x_j] = \beta_{ij}^k x_k \]
where $\beta_{ij}^k$ are the structure constants of $L$ (i.e. the matrix
elements of $\beta$).

Let $q$ be a fixed element of the ground field $\K$, $\g$ a simple Lie
algebra over $\K$, and let $\UU=\overline{U}_q(\g)$ be the
simply-connected quantised
enveloping algebra of $\g$, which is defined as follows. $\UU$ contains
a copy of the group algebra of the weight lattice of $\g$ (with basis
elements denoted by $q^\lambda$ where $\lambda$ is an integral weight of $\g$
with
respect to a Cartan subalgebra $\hh$) and generators $E_1, \ldots , E_r,
F_1, \ldots , F_r$ corresponding to fundamental roots $H_1, \ldots
, H_r$ of $\g$, with relations
\begin{eqnarray*}
q^\lambda E_i q^{-\lambda}&=&q^{\<H_i,\lambda\>}E_i\\
q^\lambda F_i q^{-\lambda}&=&q^{-\<H_i,\lambda\>}E_i\\
E_iF_j-F_jE_i &=& \delta_{ij}
   \( \frac{q^{2H_i}-q^{-2H_i}}{q_i-q_i^{-1}}\)
\end{eqnarray*}
\begin{eqnarray*}
{}[E_i,\:[E_i,\ldots ,[E_i,
E_j]_{q_i^{-n+1}}]_{q_i^{-n+3}}\cdots]_{q_i^{n-1}}&=&0\\
{}[F_i,\:[F_i,\ldots ,[F_i,
F_j]_{q_i^{-n+1}}]_{q_i^{-n+3}}\cdots]_{q_i^{n-1}}&=&0\\
\end{eqnarray*}
where $\<,\>$ denotes the Killing form in the Cartan subalgebra $\hh$
(which we do not distinguish from its dual $\hh^*$),
\[ q_i = q^{\<H_i,H_i\>},\]
\[ n = \frac{2\<H_i,H_j\>}{\<H_i,H_i\>}=\text{number
of roots $H_j + kH_i$ with $k>0$,} \]
$$ [X,Y]_p=pXY-p^{-1}YX. \leqno{\text{and}} $$
Then $\UU$  is a Hopf algebra with comultiplication
\begin{eqnarray*}
\Delta (E_i)&=&E_i\ox q^{-H_i}+q^{H_i}\ox E_i,\\
\Delta (F_i)&=&F_i\ox q^{-H_i}+q^{H_i}\ox F_i,\\
\Delta(q^\lambda )&=&q^\lambda \ox q^\lambda
\end{eqnarray*}
and antipodes
\[ S(E_i)=-q_i^{-1}E_i,\qquad S(F_i)=-q_iF_i,\qquad S(q^\lambda )=q^{-\lambda}.
\]

As in any Hopf algebra, we have the adjoint representation of $\UU$ on
itself,
given by $x\mapsto \ad x \in \End_\K \UU$ where
\[ \ad x (y)=\sum x_{(1)}yS(x_{(2)})\qquad\text{if}\quad\Delta(x)=\sum
x_{(1)}\ox
x_{(2)}. \]
This is a representation:
\begin{equation} \label{rep}
\ad (xy)=\ad x  . \ad y,
\end{equation}
and each $\ad x$ is a generalised derivation of $U$, in the sense that
\begin{equation} \label{Der1}
\ad x (yz)=\sum\ad x_{(1)}(y) . \ad x_{(2)}(z).
\end{equation}
The adjoint action of the generators of $U$ is given by
\begin{eqnarray}
\ad E_i(x)&=&E_i xq^{H_i}-q^{H_i -1}xE_i, \notag \\
\ad F_i(x)&=&F_i xq^{H_i}-q^{H_i +1}xF_i,\\
\ad q^\lambda (x)&=&q^\lambda xq^{-\lambda}. \notag
\end{eqnarray}

We use the adjoint representation to define a bracket on $\U$, defining
\begin{equation} \label{adj}
{}[x, y]=\ad x(y).
\end{equation}
Each $\ad x$ is also a generalised derivation of this bracket:
\begin{equation} \label{Jac1}
{}[x, [y, z]]=\sum [ [x_{(1)}, y], [x_{(2)}, z] ].
\end{equation}
This is a kind of Jacobi identity for the adjoint bracket. The representation
property gives another kind of Jacobi identity:
\begin{equation}\label{Jac2}
{}[\, [x, y], z]=\sum [x_{(1)}, [y, [x_{(2)}, z]\,]\,].
\end{equation}
Both of these are valid in any Hopf algebra \cite{Maj:Liealg}. However, they
are
not suitable as replacements for the Jacobi identity in defining a notion
of quantum Lie algebra, as the ad-invariant subspaces of quantised enveloping
algebras which are candidates for quantum Lie algebras are not usually
subcoalgebras, and so the coproduct which is required to formulate the
above Jacobi identities is not an intrinsic notion. That is why  we are
interested in the third kind of Jacobi identity, of the
form (\ref{Jacobi}), which is specific to quantised enveloping algebras.

The vector space $\U$ does not split up nicely into irreducible ad-invariant
subspaces; strictly speaking it is not a direct sum of such spaces, for
there are elements with components in an infinite number of them.
The {\em locally finite part} of $\UU=\U$ is the set of elements
for which this is not true, i.e.
\begin{equation}
\F ={\rm LF}\U=\{x\in \UU : \ad \UU (x) \text{ is finite-dimensional}\}
\end{equation}
The discovery of Joseph and Letzter was that $\F$ is easily described,
and has the same composition as the classical enveloping algebra $U(\g)$:
\begin{thm}[Joseph and Letzter \cite{JL}] \label{JosLet} Let $\lambda$
be a dominant integral \linebreak weight of $\g$, $V_\lambda$ the simple
$\g$-module
with highest weight $\lambda $; $V_\lambda $ can also be regarded as a
simple $\UU$-module where $\UU=\U$, and $\End_\K V_\lambda $ is a $\UU$-module
by conjugation (i.e. if $x\in \UU$ acts on $V_\lambda $ by $\rho (x)$,
it
acts on $\End _\K V_\lambda $ by $T\mapsto\sum\rho (x_{(1)})T\rho (Sx_{(2)})$).
Then
 \begin{enumerate}
  \item $E_\lambda =\ad \UU (q^{-4\lambda })$ is a finite-dimensional
$\ad$-invariant
    subspace of $\UU$, isomorphic to $\End_\K V_\lambda $ as $\UU$-module.
  \item The locally finite part of $\UU$ is
    \begin{equation}    \F=\sum_\lambda E_\lambda
    \end{equation}
    where the sum extends over all dominant integral weights of $\g$.
 \end{enumerate}
\end{thm}

For a dominant integral weight $\lambda $, we write $N_\lambda =\dim_\K
V_\lambda $. Then our result on the existence of quantum Lie algebras
in $\U$ is
\begin{thm}\label{genLiealg} For each dominant integral weight $\lambda
$, $\U$ contains an \linebreak $(N_\lambda ^2 - 1)$-dimensional $\ad$-invariant
subspace
$L_\lambda $, a linear map $\sigma :L_\lambda \ox L_\lambda \to L_\lambda
\ox L_\lambda $ and a central element $C_\lambda \in \U$ such that
\begin{eqnarray}
xy - m\circ \sigma (x \ox y) &=& C_\lambda [x, y] \qqquad (x, y\in L_\lambda
)
 \label{qcom1}\\
 \text{i.e.}\qqquad x_i x_j - \sigma _{ij}^{kl}x_kx_l &=& C_\lambda [x_i,
x_j]
   \label{qcom2}
\end{eqnarray}
where $m$ denotes multiplication in $\U$, $[\,,\,]$ is the adjoint bracket,
and   $\{x_i\}$ is a basis of $L_\lambda $.
\end{thm}
\begin{proof} Writing $\UU=\U$, let $\overline{L}=\ad \UU(q^{-4\lambda
})$. Then
$\overline{L}$ is invariant under $\ad \UU$, and by \thmref{JosLet} it
is
isomorphic
as a $\UU$-module to $\End_\K V_\lambda $. This contains a one-dimensional
submodule (spanned by the identity of $\End_\K V_\lambda $) and so by
complete reducibility of representations of $\U$ \cite[p.~324]{ChaP} we
can
write
\[  \overline{L}=\K C_\lambda \oplus L \]
as a $\UU$-module, for some $C_\lambda \in \overline{L}$. Then $C_\lambda
$ carries
the trivial representation of $\UU$:
\[  \ad x (C_\lambda )=\varepsilon (x)C_\lambda ,\qquad     \forall x
\in \UU.\]
Hence $C_\lambda $ is a central element of $\UU$, for
\[ xC_\lambda =\sum x_{(1)}C_\lambda Sx_{(2)}x_{(3)}
  =\sum\ad x_{(1)}(C_\lambda )x_{(2)}
  =\sum\varepsilon (x_{(1)})C_\lambda x_{(2)}
   =C_\lambda x .   \]
Now $\overline{L}$ is a left coideal
of $\UU$:
\begin{align}
 x \in\overline{L}\quad &\impl &\quad x &=\ad u(q^{-4\lambda })\;
   \text{ for some }u\in \UU \notag\\
 &\impl &\Delta (x)
  &=\sum\Delta (u_{(1)})\Delta (q^{-4\lambda})\Delta (Su_{(2)}) \label{Delta}\\
 &&&=\sum (u_{(1)}\ox u_{(2)})(q^{-4\lambda }\ox q^{-4\lambda })
   (Su_{(4)}\ox Su_{(3)}) \notag\\
 &&&\in \; \UU\ox \overline{L}.\notag
\end{align}
Hence for any $x\in L$ we can write
\[\Delta (x)=x_0\ox C_\lambda + \sum u^\prime \ox x^\prime \]
with $x_0, u^\prime \in U$ and $x^\prime \in L$. In fact we can choose
$C_\lambda $
so that $x_0=x$: for if $q^{-4\lambda }=C_\lambda +w$ with $w\in
L$ and if $x=\ad u(q^{-4\lambda })$, then from \eqref{Delta}
\begin{eqnarray*}
\Delta (x)&=&\sum u_{(1)}q^{-4\lambda}Su_{(3)}\ox\ad u_{(2)}(C_\lambda
+w)\\
&=&\sum u_{(1)}q^{-4\lambda} Su_{(3)}\ox\(\varepsilon(u_{(2)})C_\lambda
+ \ad u_{(2)}(w)\)\\
&=&u_{(1)}q^{-4\lambda }Su_{(2)}\ox C_\lambda +\sum u_{(1)}q^{-4\lambda
}Su_{(3)}\ox\ad u_{(2)}(w)
\end{eqnarray*}
but $\ad u_{(2)}(w)\in L$ since $L$ is invariant under $\ad U$, and the
first term is $x\ox C_\lambda $.

In $L$ we have the bracket $[x, y]=\ad x(y)$. Define $\sigma  : L\ox
L \to L\ox L$ by
\begin{equation} \label{sigma}
\sigma (x\ox y)=\sum\ad u^\prime (y)\ox x^\prime \quad
  \text{where } \Delta (x)=x\ox C_\lambda  + \sum u^\prime\ox x^\prime
\end{equation}
\begin{equation} \label{adsigma}
\text{i.e.}\qqquad\sigma (x \ox y)=\sum \ad x_{(1)}(y)\ox x_{(2)} - [x, y]\ox
C_\lambda .
\end{equation}
Then
\begin{eqnarray*}
  m\circ\sigma (x\ox y)&=&\sum x_{(1)}y(Sx_{(2)})x_{(3)} - [x,y]C_\lambda
\\[2pt]
  &=& xy - [x,y]C_\lambda .
\end{eqnarray*}
Since $C_\lambda $ is central, this establishes the theorem.
\end{proof}

The existence of a quantum Lie algebra in $U_q(\ssll (n))$ follows immediately:

\begin{thm}
If $q$ is not a $2n${\em th} root of unity nor a primitive twelfth root
of
unity, $U_q(\ssll (n))$ contains an $(n^2-1)$-dimensional quantum Lie
algebra with the adjoint bracket.
\end{thm}
\begin{proof} In Theorem 2, take $\g = \ssll (n)$ and $\lambda $ the highest
weight of the $n$-dimensional (defining) representation. Then $L_\lambda
$ carries an irreducible representation of $U_q(\g)$ (a deformation of
the adjoint representation of $\ssll (n)$). Since $C_\lambda $ is central,
$\ad C_\lambda $ acts on $L_\lambda $ as a multiple of the identity; in
a lemma we will show that this multiple is $q^2-1+q^{-2}$ (which vanishes
only if $q$ is a primitive 12th root of unity). Assuming this
for the moment, we apply ad to \eqref{qcom2} and restrict to $L=L_\lambda
$ to obtain
\begin{equation} \label{adL}
\ad_L[x_i, x_j]=\gamma _{ij}^{kl}\ad_Lx_k . \ad_Lx_l
\end{equation}
where
\[ \gamma = \frac{1-\sigma }{q^2-1+q^{-2}}. \]
Thus the quantum Jacobi identity is satisfied in $L$ by the adjoint bracket.
Also, from \ref{qcom1},
\begin{eqnarray*}
\gamma (x\ox y)=0&\impl & \sigma \(\sum x\ox y\) =\sum x\ox y\\
&\impl & C_\lambda \sum [x, y]=0\\
&\impl &\sum[x, y]=0
\end{eqnarray*}
since $\U$ contains no zero divisors \cite{JL}. Thus the adjoint bracket
restricted to $L_\lambda$ has the antisymmetry property with respect to
the quantum antisymmetriser $\gamma $.

Completion of the proof now needs the following lemmas.
\begin{lemma}\label{Clemma} The central element in $\overline{L}$ is given
by
\begin{equation}\label{C}
C_\lambda =\sum_{r=0}^{n-1}(-1)^r\frac{[n-r]_q}{[n]_q}K_r
\end{equation}
where $K_r$ is defined recursively by
\begin{equation}\label{Kr}
K_r=\ad(F_rE_r)K_{r-1},\qquad K_0=q^{-4\lambda }.
\end{equation}
\end{lemma}
\begin{lemma} \label{adc}
The adjoint action of the central element on $L$ is
\begin{equation}
\ad C_\lambda (x) = (q^2-1+q^{-2})x \quad \text{for } x\in L.
\end{equation}
\end{lemma}
The proofs, both unenlightening calculations, are relegated
to an appendix.
\end{proof}

Although the antisymmetriser $\gamma $ in this Lie algebra is not of the
Worono\-wicz form $\gamma =1-\sigma $ where $\sigma $ satisfies the braid
relation, it is of a related form. According to \ref{adL} and \ref{adsigma},
$\gamma $ is a scalar multiple of $1-\sigma $ where
\[ \sigma (x\ox y)=\overline{\sigma }(x\ox y)-[x, y]\ox C_\lambda \]
and $\overline{\sigma}
:\overline{L}\ox\overline{L}\to\overline{L}\ox\overline{L}$
is given by
\[\overline{\sigma }(x\ox y)=\sum\ad x_{(1)}(y)\ox x_{(2)}. \]
It is a straightforward exercise in Hopf algebra manipulation \cite{Wor2}
to
show that $\overline{\sigma }$ satisfies the braid relation.

\section {The quantum Lie algebra $\ssll (2)_q$}

In the case of $\g = \ssll (2)$, the simply connected quantised enveloping
algebra is generated by $E$, $F$ and $q^{\pm \half H}$ with relations
\begin{eqnarray*}
q^H Eq^{-H} &=& qE,\\
q^H F q^{-H} &=& q^{-1}F,\\
EF-FE &=& \frac{q^{2H} - q^{-2H}}{q-q^{-1}}.
\end{eqnarray*}
The highest weight of the fundamental two-dimensional representation is
(i.e. corresponds via the Killing form to) $\lambda =\half H$, so
in the notation of \lemref{Clemma} we have
\begin{eqnarray*}
K_0&=&q^{-4\lambda }\,=\,q^{-2H},\\
K_1&=&\ad(FE)K_0\,=\,-(q-q^{-1})(qEF - q^{-1}FE).
\end{eqnarray*}
Hence, according to \lemref{Clemma}, the central element is
\begin{equation}
C=q^{-2H}+\(\frac{q-q^{-1}}{q+q^{-1}}\) (qEF-q^{-1}FE)
\end{equation}
A basis for the quantum Lie algebra is $(X_\pm ,\:X_0)$ where
\begin{eqnarray}
X_+&=& \frac {\ad E(q^{-2H})}{q-q^{-1}}=q^{-H}E,\notag\\
X_- &=& -\frac{\ad F(q^{-2H})}{q-q^{-1}}=q^{-H}F,\\
X_0 &=& - \frac{\ad (FE)q^{-2H}}{q^2-q^{-2}}=\frac{qEF-q^{-1}FE}{q+q^{-1}}
   \notag\\
&=&\frac{C-q^{-2H}}{q-q^{-1}}. \notag
\end{eqnarray}
The adjoint brackets of these basis elements can be calculated to be
\begin{equation}\label{[]}
  \begin{array}{ccc}
    {[}X_{+},X_{+}]=0,     &  [X_{+},X_{0}]=-q^{-1}X_{+},   &
      [X_{+},X_{-}]=(q+q^{-1})X_{0}, \\[5pt]
    {[}X_{0},X_{+}]=qX_{+},& [X_{0},X_{0}]=(q-q^{-1})X_{0},&
      [X_{0},X_{-}]=-q^{-1}X_{-}, \\[5pt]
    {[}X_{-},X_{+}]=-(q+q^{-1})X_{0}, & \:[X_{-},X_{0}]=qX_{-}, &
      [X_{-}, X_{-}]=0.
  \end{array}
\end{equation}
The quantum Lie algebra $\ssll (2)_q$ is defined by these brackets together
with the quantum antisymmetriser $\gamma =(q^2-1+q^{-2})^{-1}\gamma ^\prime$
where
\begin{equation}\label{qantisym}
  \begin{split}
   \gamma ^\prime (X_\pm\ox X_\pm)&=0,\\[3pt]
   \gamma ^\prime (X_\pm\ox X_0)&=q^{\mp 2}X_\pm\ox X_0
    - X_0\ox X_\pm,\\[3pt]
   \gamma ^\prime (X_0\ox X_\pm)&=q^{\pm 2}X_0\ox X_\pm
    - X_\pm\ox  X_0,\\[3pt]
   \gamma ^\prime (X_+\ox X_-)&=X_+\ox X_- - X_ -\ox X_+
    +(q^2-q^{-2})X_0\ox X_0,\\
   \gamma ^\prime (X_-\ox X_+)&=-\gamma^\prime(X_+\ox X_-),\\
   \gamma ^\prime (X_+\ox X_0)&=(q-q^{-1})^2X_0\ox X_0
    +\(\frac{q-q^{-1}}{q+q^{-1}}\) (X_+\ox X_- - X_-\ox X_+)\\
   &=\(\frac{q-q^{-1}}{q+q^{-1}}\)\gamma ^\prime (X_+\ox X_-).
  \end{split}
\end{equation}
We note the following properties of $\ssll (2)_q$:
\begin{propn}
 \begin{enumerate}
  \item $\ssll (2)_q$ is a balanced quantum Lie algebra.
  \item The quantum antisymmetriser $\gamma$ is essentially idempotent:
   \[ \gamma ^2 = \(\frac{q^2+q^{-2}}{q^2-1+q^{-2}}\)\gamma . \]
 \end{enumerate}
\end{propn}
\begin{proof} By calculation. \end{proof}
The second of these properties is peculiar to $\ssll (2)_q$, reflecting
the simple structure of its representations, as will become apparent when
we consider $\ssll (3)_q$. The first, however, may be more general.

The enveloping algebra of the quantum Lie algebra $U(\ssll (2)_q)$ --- not
to be confused with the quantised enveloping algebra $U_q(\ssll (2))$ ---
is defined by means of the brackets \eqref{[]} and the antisymmetriser
\eqref{qantisym}.
By redefining the generators to absorb the factor $q^2-1+q^{-2}$, it can
be presented as the associative algebra generated by three elements
$Y_0,\;Y_\pm$
with relations
\begin{equation}\label{Y}
  \begin{split}
    qY_0Y_+ - q^{-1}Y_+Y_0 &=Y_+\\[3pt]
    q^{-1}Y_0Y_- - qY_-Y_0&=-Y_-\\[3pt]
    Y_+Y_--Y_-Y_++(q^2-q^{-2})Y_0^2&=(q+q^{-1})Y_0.
  \end{split}
\end{equation}

On the other hand, as elements of the quantised enveloping algebra $U_q(\ssll
(2))$ the generators $X_\pm$, $X_0$ and the central element $C$ satisfy
(by \linebreak \thmref{genLiealg})
\begin{eqnarray}\label{XC}
  q^2X_0X_+-X_+X_0&=&qCX_+,\notag\\[3pt]
  q^{-2}X_0X_--X_-X_0&=& - q^{-1}CX_-,\notag\\[3pt]
  X_+X_--X_-X_++(q^2-q^{-2})X_0^2&=&(q+q^{-1})CX_0\\[3pt]
\text{and}\qqqquad  CX_m=X_mC\qquad &\phantom{=}& \qquad (m=0,\: \pm 1).\notag
\end{eqnarray}
However, the central element $C$ is not independent of the other generators
$X_m$:
there is a further relation
\begin{equation}\label{Cas}
C^2=(q-q^{-1})^2\( X_0^2 + \frac{qX_-X_++q^{-1}X_+X_-}{q+q^{-1}}\) +1
\end{equation}
as can be verified by direct calculation.

The main results of this section are that the locally finite part of $\U$
is isomorphic, as an algebra, to the algebra with generators $X_0,\: X_\pm
,\: C$ and relations \eqref{XC}, and as a vector space to the algebra of
polynomials in four commuting variables $X_0,\: X_\pm,\: C$ subject to
the relation \eqref{Cas}. We use the following notation for the various
algebras:

$\A$ = enveloping algebra of the quantum Lie algebra $\ssll (2)_q$

{}\quad = algebra generated by $Y_\pm ,\; Y_0$ with relations \eqref{Y},

$\B$ = algebra generated by $X_\pm,\; X_0,\; C$ with relations \eqref{XC},

$\C$ = algebra generated by $X_\pm,\; X_0,\; C$ with relations \eqref{XC}
and \eqref{Cas},

$\F$ = locally finite part of $\U$.
\newline
Then $\A$, $\B$ and $\C$ are related by homomorphisms
\begin{eqnarray*}
 \varphi &:& \B \to \A\quad \text{with  } \varphi (C) =1\\
 \psi &:& \B\to\C\quad \text{enforcing \eqref{Cas}}.
\end{eqnarray*}
\begin{lemma} \label{lem1} The algebra $\A$ has a basis of ordered monomials
$Y_-^lY_0^mY_+^n$ where $l,\: m,\: n$ are non-negative integers. \end{lemma}
\begin{proof} It is evident that the commutation relations \eqref{Y} allow
any product of $Y$s to be expressed as a combination of ordered monomials.
We have only to show that these monomials are independent. According to
the diamond lemma \cite{Berg}, this will follow if it can be shown that
the two methods of reducing $Y_+Y_0Y_-$ to a combination of ordered monomials
both lead to the same result. This is easily checked using \eqref{Y}.
\end{proof}

\begin{lemma}\label{lem2} The algebra $\B$ has a basis of ordered monomials
$C^kX_-^lX_0^mX_+^n$. \end{lemma}
\begin{proof} $\A\ox\K[C]$ contains an isomorphic copy of $\B$ generated
by $C$ and $X_m=CY_m$. This is spanned by monomials
$C^kX_-^lX_0^mX_+^n=C^{k+l+m+n}Y_-^lY_0^mY_+^n$
and by \lemref{lem1} these are independent. \end{proof}

\begin{lemma}\label{lem3} The algebra $\C$ is isomorphic as a vector space
to the space of polynomials in four commuting variables $C,\; X_0,\; X_+,\;
X_-$ subject to the relation \eqref{Cas}. \end{lemma}
\begin{proof} Let
\[   C_2=(q+q^{-1})X_0^2+qX_-X_++q^{-1}X_+X_-.  \]
Using \eqref{XC}, it can be checked that $C_2$ commutes with $X_0,\; X_+$
and $X_-$. Hence the two-sided ideal $I$ generated by
$(q+q^{-1})(C^2-1)-(q-q^{-1})^2C_2$
(i.e. the kernel of the relation \eqref{Cas}) is the same as the
left ideal generated by it. It follows that the monomials $X_-^lX_0^mX_+^n$
and $CX_-^lX_0^mX_+^n$ are independent modulo $I$, and therefore form
a basis of $\C$. But these monomials also form a basis of the commutative
algebra described in the theorem. \end{proof}

\begin{thm} \label{locfin}
The locally finite part of $\Utwo$ is isomorphic to $\C$. \end{thm}
\begin{proof} By \thmref{JosLet}, the locally finite part of $\UU=\U$
is
\[ \F = \sum_{n=0}^\infty V_n \qquad \text{where}\quad V_n=\ad \UU(q^{-2nH}).
\]
Since $\half nH$ is the highest weight of the $(n+1)$-dimensional
representation of $\UU$, $\dim V_n=(n+1)^2$. By the derivation property
\eqref{Der1} of the adjoint action,
\[ V_{m+n} \subseteq V_mV_n. \]
$$  V_n \subseteq V_1^n. \leqno{\text{Hence}} $$
Let $\UU^\prime$ be the subalgebra of $\UU$ generated by $V_1$, i.e. by
$C,\:
X_-,\: X_+$ and $X_0$. Since the relations \eqref{XC} and \eqref{Cas}
hold
in $\UU$, $\UU^\prime$ is a quotient of $\C$ and therefore, by \lemref{lem3},
is spanned by monomials $X_-^kX_0^lX_+^m$ and $CX_-^kX_0^lX_+^m$. Let
$W_n$ be the subspace spanned by monomials $X_-^kX_0^lX_+^m$ with
$k+l+m=n$. Then
\begin{equation}\label{V1}
V_1^{n+2}=V_1^n+CW_{n+1} +W_{n+2}
\end{equation}
as can be proved by induction, using
$$  V_1W_r =CW_r +W_{r+1}  $$
$$ C^2 W_r +W_{r+2}=W_r + W_{r+2}   \leqno{\text{and}} $$
which follows from \eqref{Cas}. But
\begin{equation}\label{V3}
 \dim W_r\leq \half (r+1)(r+2),
\end{equation}
$$ \dim(V_1^{n+2})\leq\dim (V_1^n)+(n+3)^2.  \leqno{\text{so}} $$
We can now prove by induction that
\[ V_1^n = V_n\oplus V_{n-2}\oplus \cdots \oplus (V_0 \text{ or } V_1).
\]
For if this is true, then from \eqref{V1} we have
\[ V_1^{n+2}\supseteq V_{n+2}\oplus V_n \oplus \cdots \]
since the $V_n$ are known to be complementary in $\UU$ \cite{JL}. By dimensions
the inclusion must be equality, so the sum in \eqref{V1} is direct and
also
\eqref{V3} must be equality. It follows that the monomials $X_-^kX_0^lX_+^m$
and $CX_-^kX_0^lX_+^m$ are independent and form a basis of $\sum V_n=\F$.
By \lemref{lem3}, $\F$ is isomorphic to $\C$. \end{proof}

In one sense, the locally finite part of the quantised enveloping algebra
is the same size as the classical enveloping algebra of $\ssll (2)$.
Both can be written as $Z\ox \UU_0$ where $Z$, the centre, is a polynomial
algebra in one variable and $\UU_0$ has the same decomposition into
simple modules under the adjoint action in both cases. However, the above
shows that it is more appropriate to regard the locally finite part of
$\Utwo$ as twice as big as the classical enveloping algebra. In the classical
case the centre is generated by a Casimir element which is a polynomial
in the Lie algebra generators; in the quantum case the centre is generated
by a square root of such a polynomial, tending to 1 in the classical limit.

\section{The Quantum Lie Algebra $\ssll (3)_q$}

This section contains some details of and observations on the quantum
Lie algebra of type $A_2$. All were established by mindless calculation.

We take $\lambda $ to be the highest weight of the triplet representation,
\[ \lambda =\third (2H_1 + H_2) \]
(identifying weights with elements of the Cartan subalgebra by means of
the Killing form, and taking $H_1$ and $H_2$ to have unit length). A convenient
basis for the quantum Lie algebra $\ssll (3)_q$ is the following, consisting
of two elements $T_1, T_2$ (scalar multiples of the $K_1, K_2$ defined
in \lemref{Clemma}) and elements $X_\alpha = X_{\pm 1},\,X_{\pm2},\,X_{\pm 12}$
corresponding to the roots $\alpha = \pm H_1,\,\pm H_2,\,\pm (H_1+H_2)$ of
$\ssll (3)$:
\begin{eqnarray*}
T_1&=&\frac{\ad (F_1E_1)q^{-4\lambda }}{q-q^{-1}}
    = q^{-\thrd{2}(H_1+2H_2)} (qE_1F_1-q^{-1}F_1E_1)\\
T_2&=&\ad (F_2E_2)T_1\\
X_1&=&\frac{\ad E_1(q^{-4\lambda})}{q-q^{-1}}
    =q^{-\third(5H_1+4H_2)}E_1 \\
X_{12}&=&\ad E_2(X_1)
    =q^{-\third(5H_1+H_2)}(E_2E_1-q^{-1}E_1E_2)\\
X_2&=&\ad F_1(X_{12})\\
%   &=&q^{-\third(2H_1+H_2)}
%      (-qE_2E_1F_1 +E_1E_2F_1 +q^{-1}F_1E_2E_1 -q^{-2}F_1E_1E_2)\\
X_{-1}&=&\frac{\ad F_1(q^{-4\lambda})}{q-q^{-1}}
       =-q^{-\third(5H_1+4H_2)}F_1\\
X_{-12}&=&\ad F_2(X_{-1})
    =q^{-\third(5H_1+H_2)}(qF_1F_2-F_2F_1)\\
X_{-2}&=&-\ad E_1(X_{-12})\\
%     &=&q^{-\third(2H_1+H_2)}
%      \left\{q(qE_1F_1-q^{-1}F_1E_1)F_2-F_2(qE_1F_1-q^{-1}F_1E_1)\right\}
\end{eqnarray*}
The adjoint action of $\U$ on these basis elements is:
\begin{eqnarray*}
   \ad(q^{H_i})T_j&=&T_j\\
   \ad(q^{H_i})X_{\pm j}&=&\begin{cases}
     q^{\pm 1}X_{\pm j}   &\text{if $i=j$}\\
     q^{\mp 1/2 }X_{\pm j} &\text{if $i\ne j$}
       \end{cases}\\
  \ad(q^{H_i})X_{\pm 12}&=&q^{\pm 1/2}X_{\pm 12}
\end{eqnarray*}
\begin{equation*}
  \ad E_i(X_\alpha )=\begin{cases}
      X_{\alpha +H_i}  &\text{if $\alpha +H_i$ is a root}\\
      T_i &\text{if $\alpha =-H_i$ }\\
       0  &\text{otherwise} \end{cases}
\end{equation*}
\begin{equation*}
 \ad F_i(X_\alpha )=\begin{cases}
    X_{\alpha -H_i}  &\text{if $\alpha -H_i$ is a root}\\
    T_i &\text{if $\alpha =H_i$}\\
    0  &\text{otherwise} \end{cases}
\end{equation*}
\begin{equation*}
  \ad E_i(T_j)=\begin{cases}
    (q+q^{-1})X_i  &\text{if $i=j$}\\
    X_{12}  &\text{if $i\ne j$} \end{cases}
\end{equation*}
\begin{equation*}
 \ad F_i(T_j)=\begin{cases}
    (q+q^{-1})X_{-i}  &\text{if $i=j$}\\
    X_{-12}  &\text{if $i\ne j$} \end{cases}
\end{equation*}
{}From this the adjoint brackets between the basis elements can be calculated;
the result is given in the table, in which $[X,Y]$ is given in the row
labelled
by $X$ and the column labelled by
$Y$. If $\alpha +\beta \ne 0$, $[X_\alpha , X_\beta ]=c_{\alpha \beta
}X_{\alpha +\beta}$ is represented by the coefficient $c_{\alpha \beta }$;
similarly, $[T_i, X_\alpha ] = b_{i\alpha }X_\alpha $ is represented
by $b_{i\alpha }$.
\begin{figure}[h]
\begin{tabular}{c||c|c|c|c|}
  &$T_1$ & $T_2$ & $X_1$ & $X_{-1}$\\ \hline \hline
  $T_1$ & $-(q^2-q^{-2})T_1$ & $-(q-q^{-1})T_1$ & $-q(q+q^{-1})$
    & $q^{-1}(q+q^{-1})$\\ \hline
  $T_2$ & $-(q-q^{-1})T_1$ & & $-q$ & $q^{-1}$\\ \hline
  $X_1$ & $q^{-1}(q+q^{-1})$ & $q^{-1}$ & 0 & $T_1$\\ \hline
  $X_{-1}$ & $-q(q+q^{-1})$ & $-q$ & $-T_1$ & 0\\ \hline
  $X_2$ & $-q^2$ & $-q^{-1}(q^2+q^{-2})$ & $-q^{3/2}$ & 0\\ \hline

  $X_{-2}$ & $q^{-2}$ & $q(q^2+q^{-2})$ & 0 & $q^{-3/2}$\\ \hline

  $X_{12}$ & 1 & $-q^{-3}$ & 0 & $q^{1/2}$\\ \hline
  $X_{-12}$ & $-1$ & $q^3$ & $-q^{-1/2}$ & 0\\ \hline
\end{tabular}
\noindent\par\medskip
\begin{tabular}{c||c|c|c|c|}
 &$X_2$ & $X_{-2}$ & $X_{12}$ & $X_{-12}$\\ \hline \hline
$T_1$ & $q^{-2}$ & $-q^2$ & $-1$ & 1\\ \hline
$T_2$ & $q(q^2+q^{-2})$
    & $-q^{-1}(q^2+q^{-2})$ & $q^3$ & $-q^{-3}$\\ \hline
$X_1$ & $q^{-3/2}$ & 0 & 0 & $q^{1/2}$\\ \hline
$X_{-1}$ & 0 & $-q^{3/2}$ & $-q^{-1/2}$ & 0\\ \hline
$X_2$ & 0 & & 0 & $-q^{-5/2}$\\ \hline
$X_{-2}$ & & 0 & $q^{5/2}$ & 0\\ \hline
$X_{ 12}$ & 0 & $-q^{-5/2}$ & 0 & $-q^{-1}T_1+T_2$\\ \hline
$X_{-12}$ & $q^{5/2}$ & 0 & $qT_1-T_2$ & 0\\ \hline
\end{tabular}
\caption{The brackets of $\ssll (3)_q$.}
\end{figure}

The entries which are too long to fit into the table are
\begin{eqnarray*}
  {[}T_2, T_2] &=& -(q^2-q^{-2})T_1+(q^3-q^{-3})T_2,\\
  {[}X_2, X_{-2}] &=& q^{-1}(q-q^{-1})T_1 - qT_2,\\
  {[}X_{-2}, X_2] &=& q(q-q^{-1})T_1+q^{-1}T_2.
\end{eqnarray*}
Notice that if $L$ is regarded as a module over $\K [q, q^{-1}]$ then, with
this basis, the brackets exhibit the simple antisymmetry which was also
found by Delius and H\"uffmann. Define a map of {\em $q$-conjugation} on
$L$, $X \to X^\blacktriangledown$, by the requirements that it leaves
the above
basis elements fixed and that
$\{f(q)X\}^\blacktriangledown=f(q^{-1})X^\blacktriangledown$.
Then it is apparent from the table that
\begin{equation}
[Y, X]=-[X, Y]^\blacktriangledown \quad \text{for all } X, Y \in L.
\end{equation}

  The quantum antisymmetriser of $\U$ can be described as follows. The
quantum Lie algebra $L=\ssll (3)_q$ carries an irreducible representation
of $\U$; as in the classical case, $L \ox L$ splits into irreducible
representations
as
\[ L \ox L = M_{27} \oplus M_{10} \oplus M_{10}^* \oplus
 M_8^{\text{a}} \oplus M_8^{\text{s}} \oplus M_1 \]
where the numerical labels 27, 10, 8, 1 indicate the dimensions of the
simple modules. These are completely defined by their highest-weight elements:
\begin{equation}\label{hweight}
 \begin{split}
W_{27}&=X_{12}\ox X_{12},\\
W_{10}&=q^{1/2}X_1\ox X_{12} - q^{_1/2}X_{12}\ox X_1,\\
W_{10}^* &= q^{1/2}X_2 \ox X_{12} - q^{-1/2}X_{12}\ox X_2,\\
W_8^{\text{s}} &= \left\{q^{5/2}(q+q^{-1})+q^{-5/2}\right\}X_1 \ox X_2
+ \left\{q^{-5/2}(q+q^{-1})+q^{5/2}\right\}X_2
\ox X_1\\
        &\quad -(q^4T_1+q^{-1}T_2)\ox X_2 -X_2 \ox (q^{-4}T_1+qT_2),\\
W_8^{\text{a}} &= q^{3/2}X_1\ox X_2 - q^{-3/2}X_2 \ox X_1
  -q(qT_1-T_2)\ox X_{12}\\
  &\quad + q^{-1}X_{12}\ox (q^{-1}T_1 - T_2),\\
W_1&=T_1\ox T_2 +T_2\ox T_1 -(q+q^{-1})(T_1\ox T_1 +T_2\ox T_2)\\
  &\quad +(q^2+1+q^{-2})(q^{-1}X_1\ox X_{-1} +qX_{-1}\ox X_1\\
   &\qquad +q^{-1}X_2\ox X_{-2} + qX_{-2}\ox X_2
    -q^{-2}X_{12}\ox X_{-12} - q^2 X_{-12}\ox X_{12}).
\end{split}
\end{equation}

  Let $\Pi _{27}, \Pi _{10}, \Pi _{10}^*, \Pi _8^{\text{s}},
\Pi
_8^{\text{a}}$ and $\Pi _1$ denote the projectors onto these simple modules.
Then
\begin{lemma} The quantum antisymmetriser of $\U$ is
\begin{equation}
\gamma = \frac
  {(q^2+1)\Pi _{10}^* + (q^2+q^{-2})\Pi _8^{\text{a}} +(1+q^{-2})\Pi
_{10}}
  {q^2 - 1 + q^{-2}}
\end{equation}
\end{lemma}
\begin{proof} Straightforward Hopf-algebra manipulations show that $\gamma
$ commutes with the diagonal adjoint action of $\mathcal{U}=U_q(\ssll
(3))$ on $L \ox L$, i.e. with $(\ad \ox \ad )\circ \Delta (u)$ for all
$u \in \mathcal{U}$. Thus it is only necessary to calculate the action
of $\gamma $ on the highest-weight elements listed above. This is not
too laborious for the cases of $W_{10}, W_{10}^*$ and $W_{27}$.
For the octets the relevant highest-weight elements (eigenvectors of $\gamma
$) and the eigenvalues can be determined from the Jacobi identity: in
the two-dimensional space of highest-weight octet vectors (elements of
$L\ox L$ with weight $H _1 + H _2$ which are annihilated by
$(\ad\ox\ad)\circ\Delta (E_1)$ and $(\ad\ox\ad)\circ\Delta (E_2)$, the
eigenvectors of $\gamma $ are the relative eigenvectors of $\ad\circ\beta $ and
$m\circ(\ad\ox\ad)$.

The eigenvalue 0 for the singlet can be verified
in a similar way. We find that $\beta (W_1)=0$ but $m\circ (\ad\ox\ad
)(W_1)\ne 0$. Since $\ad\circ\beta (W_1)=m\circ (\ad\ox\ad )\circ\gamma
(W_1)$ and $\gamma (W_1)$ must be a multiple of $W_1$, it follows that
$\gamma (W_1)=0$. \end{proof}

The quantised enveloping algebra $\UU=\Uthree$ has a symmetry between
$(E_1,F_1,H_1)$
and $(E_2,F_2,H_2)$ corresponding to the outer automorphism of $\ssll
(3)$ which comes from the symmetry of the Dynkin diagram $\circ$---$\circ$.
The quantum Lie algebra $\L=\ssll (3)_q=\ad\UU(q^{-4\lambda })$ is not
invariant under this automorphism of $\UU$. Indeed, the weight diagram
of the 3-dimensional representation is taken by this symmetry to that
of the conjugate representaion, with highest weight $\lambda ^*
\ne\lambda $. This gives another quantum Lie algebra
$L^*\subset\ad\UU(q^{-4\lambda ^*})$
with $L^*\ne L$. This is another octet (under the adjoint representation)
which also generates the locally finite part of $\UU$. It is interesting to
express the elements
of $L^*$ in terms of those of $L$ and to compare this situation with the
classical one.

In the classical enveloping algebra the octets are the Lie algebra itself,
one which is quadratic in the generators, and multiples of these by functions
of the Casimirs.
The outer automorphism acts linearly on the generators. In the locally
finite part of the quantised
enveloping algebra we have the same general structure \cite{JL}; the quadratic
octet has highest-weight element obtained by multiplying the factors of
the tensor $W_8^{\text{s}}$ in \eqref{hweight}. The highest-weight element
of the octet $L$ is $X_{12}$. Let $X_{12}^*$ be the highest-weight
element of $L^*$, and let $Y_{12}=m(W_8^{\text{s}})$. Then a
calculation
gives
\begin{equation}
X_{12}^*=\frac{q^{1/2}(q+q^{-1})CX_{12}+q(q-q^{-1})Y_{12}}
                            {(q^{1/2}+q^{-1/2})(q^2+1+q^{-2})}.
\end{equation}
Thus the elements of $L^*$ are quadratic in the generators $L$,
and the quantum counterpart of the outer automorphism is a nonlinear symmetry
of $\Uthree$.

\section{Appendix. Proofs of Lemmas \ref{Clemma} and \ref{adc}}
\noindent \textbf {Lemma \ref{Clemma}} The central element in $\overline{L}$
is
given by
\begin{equation}\tag{\ref{C}}
C_\lambda =\sum_{r=0}^{n-1}(-1)^r\frac{[n-r]_q}{[n]_q}K_r
\end{equation}
where $K_r$ is defined recursively by
\begin{equation}\tag{\ref{Kr}}
K_r=\ad(F_rE_r)K_{r-1},\qquad K_0=q^{-4\lambda }.
\end{equation}

\begin{proof}  Let $C$ be defined by \eqref{C}. Since this has zero weight
for the Cartan algebra generated by $H_1,\ldots ,H_{n-1}$ (i.e. $\ad
q^{H_i}(C)=C$),
in order to show that it is central it is sufficient to show that $\ad
E_i(C)=0$ for $i=1,\ldots ,n-1$. We calculate $\ad E_i(K_j)$.

  Since $\ad \U (q^{-4\lambda })$ is the direct sum of the $q$-analogues
of the scalar and the adjoint modules of $\ssll (n)$, the structure of
the roots of $\ssll (n)$ gives
\begin{equation}\label{adE}
\ad E_j \ad E_i \ad E_{i-1} \cdots \ad E_1(q^{-4\lambda })=0
  \quad \text{unless }j=i+1.
\end{equation}
$$  \ad E_j(K_i)=0 \quad \text{if } j > i+1 \text{ and } i > 0.
  \leqno{\text{Hence}} $$
This also holds if $i=0$ since $q^{-4\lambda }$ commutes with $E_j$ for
$j > 1$.

  For $i=j$ we have
\begin{eqnarray*}
\ad E_i(K_i) &=& \ad \( F_iE_i + [2H_i]_q \)\ad E_i(K_{i-1}\\
  &=& [\< 2H_i, H_i\>]_q\ad E_i(K_{i-1})\\
  &=& [2]_q\ad E_i(K_{i-1})
\end{eqnarray*}
since
\begin{eqnarray*}
  (\ad E_i)^2K_{i-1} &=& \ad F_{i-1} \ldots \ad F_1 (\ad E_i)^2 \ad E_{i-1}
   \ldots \ad E_1(q^{-4\lambda })\\
  &=& 0 \qquad \text{by \eqref{adE}}.
\end{eqnarray*}

  For $j = i - 1$ we have
\begin{eqnarray} \label{FEEK}
\ad E_j(K_{j+1}) &=&
    \ad F_{j+1}\ad ( F_jE_j + [2H_j])\ad E_{j+1}\ad E_j(K_{j-1})\notag\\
  &=& \ad ( [2H_j - 1]_q)\ad F_{j+1}\ad E_{j+1} \ad E_j(K_{j-1})\notag\\
  &=& \ad (F_{j+1}E_{j+1}E_j)K_{j-1}
    \qqquad\text{since } [2\< H_j, H_j\> - 1]_q =1\notag\\
  &=& \ad [(E_{j+1}F_{j+1} - [2H_{j+1}]_q)
     E_j] K_{j-1}\notag\\
  &=& \ad E_j(K_{j-1})
\end{eqnarray}
since $ \ad F_{j+1} \ad E_j(K_{j-1}) = 0$
by the structure of the roots of $\ssll (n)$.

\smallskip

  For $j < i - 1$,
\begin{eqnarray*}
 \ad E_j(K_i) &=& \ad (F_i \ldots F_{j+1}E_j F_j \ldots F_1 E_i \ldots
E_1)
   q^{-4\lambda }\\[3pt]
 &=& \ad [ F_i\ldots F_{j+1}\( F_jE_j + [2H_j]_q\)
   F_{j-1}\ldots F_1 E_i\ldots E_1] q^{4\lambda }\\[3pt]
 &=& [\< 2H_j, H_{j+1}+\cdots+H_i\> ]_q
   \ad\left( F_i E_i\ldots F_{j+1}E_{j+1}E_j\right) K_{j-1}\\
 &&\qqqquad\qqqquad\qqqquad\qqqquad\qquad\text{ by \eqref{adE}}\\[3pt]
&=& - \ad \left( F_iE_i\ldots F_{j+2}E_{j+2}E_j\right) K_{j-1}
    \qqqquad\text{ by \eqref{FEEK}}\\[3pt]
&=& 0 \qqqquad \text{ by \eqref{adE}.}
\end{eqnarray*}

  To summarise,
\begin{equation}
 \begin{split}
  \ad E_j(K_i) &= 0 \qqquad \text{if $i < j-1$ or $i > j+1$,}\\
  \ad E_j(K_{j+1}) &= \ad E_j (K_{j-1})\\
  \text{and }\qqqquad \ad E_j (K_j) &= [2]_q \ad E_j (K_{j - 1}).
 \end{split}
\end{equation}

  Hence for $j=1,\ldots, n - 2$,
\begin{eqnarray*}
\ad E_j (C) &=& \frac{(-1)^{j-1}}{[n]_q}
    \( [n-j-1]_q - [2]_q[n - j]_q + [n-j+1]_q \)\ad E_j (K_{j-1})\\
& =& 0,
\end{eqnarray*}
while
\[ \ad E_{n-1}(C) = \frac{(-1)^{n-1}}{[n]_q}\([2]_q - [2]_q[1]_q \) =
0. \]
Thus $C$ is a highest-weight vector in the module $\overline{L}$. Since
its weight is zero, it follows that $C$ is an invariant element ($\ad
x(C)=\varepsilon (x)C$ for all $x\in\UU$) and therefore central. It is
therefore a multiple of $C_\lambda$. But $C=q^{-4\lambda}+w$ with $w\in
L$; hence $C=C_{\lambda}$.
\end{proof}

\noindent \textbf{\lemref{adc}} The adjoint action of the central element
on
$L$ is
\begin{equation}
\ad C_\lambda (x) = (q^2-1+q^{-2})x \quad \text{for } x\in L.
\end{equation}

\begin{proof} It follows from \lemref{Clemma} by Schur's lemma that $\ad
C_\lambda
$ acts as a multiple of the identity on the irreducible component $L$.
To find the multiple we calculate $\ad C_\lambda (X_1)$ for $X_1=\ad
E_1(q^{-4\lambda
})$. We use
\begin{eqnarray*}
\ad E_r(X_1) &=& 0 \quad \text{unless  } r=2,\\
\ad F_r(X_1) &=& 0 \quad \text{unless  } r=1,\\
\ad F_1(X_1) &=& K_1.
\end{eqnarray*}
Then we have
\[ \ad K_0(X_1) = q^{-\< 4\lambda , H_1\>}X_1 = q^{-2}X_1. \]
Now
\[ \ad E_1 (q^{-4\lambda }) =  (q-q^{-1})q^{-4\lambda + H_1}E_1, \]
so \begin{eqnarray}\label{K1}
 K_1 &=& (q-q^{-1})\ad F_1(q^{-4\lambda +H_1}E_1)\nonumber\\
  &=& - (q-q^{-1})q^{-4\lambda +2H_1}(qE_1F_1 - q^{-1}F_1E_1)\nonumber\\
  &=&(q-q^{-1})^2q^{-4\lambda +2H_1}F_1E_1 - q.q^{- 4\lambda }(q^{4H_1
- 1})
\end{eqnarray}
Therefore
\begin{eqnarray}\label{K1X1}
\ad K_1 (X_1) &=& - q\( q^{\<-4\lambda +4H_1, H_1\>}
   - q^{\<-4\lambda , H_1\>}\)X_1 \nonumber \\
&=& -q(q^2 - q^{-2})X_1.
\end{eqnarray}

  To calculate $\ad K_2(X_1)$, note that
\[ \ad E_2 (K_1) = E_2 K_1 q^{H_2} - q^{H_2 - 1}K_1E_2
  = q^{H_2 - 1}(E_2 K_1 - K_1 E_2). \]
Since $\ad F_2 (X_1) = 0$, we have
\begin{eqnarray*}
  \ad K_2 (X_1) &=& \ad \left[ F_2.\ad E_2(K_1).q^{H_2}\right] X_1\\
   &=& \ad \left[ q^{2H_2}F_2(E_2K_1 - K_1E_2)\right] X_1.
\end{eqnarray*}
Using \eqref{K1X1} and the fact that $\ad F_2(X_1)=0$,
\begin{eqnarray*}
  \ad\left[ q^{2H_2}F_2K_1E_2\right] X_1 &=& - q\ad\left[ q^{2H_2}
   \( E_2F_2 - [2H_2]_q\) \right] X_1\\
  &=& - (q^2 - q^{-2})X_1;
\end{eqnarray*}
using \eqref{K1} and the fact that $\ad(E_1E_2)X_1=0$,
\begin{eqnarray*}
  \ad\left[ q^{2H_2}F_2K_1E_2\right] X_1 &=&
   - q\ad\left[ q^{2H_2}F_2q^{-4\lambda }(q^{4H_1} - 1)E_2\right]X_1\\
  &=& - q\ad\left[ q^{2H_2 - 4\lambda }\( q^{4H_1 - 2} - 1\)
    \( E_2F_2 - [2H_2]_q\) \right]X_1\\
  &=& - q^{-2}(q^2 - 1)X_1;
\end{eqnarray*}
\[ \therefore \qquad \ad K_2(X_1) = - (q^2 - 1)X_1. \]

  We now have that $K_3$ is a sum of products of $\ad E_3$ and $\ad F_3$,
both of which annihilate $X_1$, and operators of which $X_1$ is an eigenvector;
it follows that $\ad K_3(X_1)=0$. Similarly $\ad K_r(X_1)=0$ for $r > 3$; hence
\begin{eqnarray*}
\ad C_\lambda (X_1) &=& \(q^{-2} + q(q-q^{-1})\frac{(q+q^{-1})[n-1]_q -
[n-2]_q}{[n]_q}
   \)X_1\\
  &=&(q^2 - 1 + q^{-2})X_1.
\end{eqnarray*}
\end{proof}

\bibliographystyle{plain}

\end{document}